\documentclass[xypic,amscd,syntonly,amssymb,verbatim,12pt]{amsart}

\usepackage{graphicx}
\usepackage[all]{xy}
\usepackage{amssymb}
\usepackage{amsmath}
\usepackage[mathscr]{euscript}

\usepackage{epsfig}
\theoremstyle{plain}
\newtheorem{Thm}{Theorem}
\newtheorem{Prop}[Thm]{Proposition}

\newtheorem{Def}[Thm]{Definition}

 \theoremstyle{definition}
\theoremstyle{remark}
\newtheorem{Rems} {Remarks}
\numberwithin{equation}{section}

\begin{document}
 %\title{Coupling class of actions of reductive groups}
 %\title{Vertex operators between general $B$-branes}
  \title{ Holomorphic gauge fields on $B$-branes}

 \author{ ANDR\'{E}S   VI\~{N}A}
\address{Departamento de F\'{i}sica. Universidad de Oviedo.   Avda Calvo
 Sotelo 18.     33007 Oviedo. Spain. }
 \email{vina@uniovi.es}
 
%\thanks{This work has been partially supported by Ministerio de Ciencia y
%Tecnolog\'{\i}a, grant FPA2009-11061}

 \maketitle
 
 \centerline{ In memoriam to Ange Vi\~{n}a }
 
\begin{abstract}
%%%%%%%%%%%%%%%%%%%%%%%%%%%%%%%%%%%%%%%%%%%%%%%%%%%%%%%%

 Considering the $B$-branes over a complex manifold   as the objects of 
%the corresponding
 the bounded derived category of coherent sheaves on that manifold,  we extend the definition of holomorphic gauge fields on vector bundles   to $B$-branes.
We construct a family   of coherent sheaves on the complex projective space, 
%${\mathbb P}^n$,  
 which generates the corresponding bounded derived category 
 %on the projective space
  and  such that
the supports of the elements of this family are two by two disjoint.
Using that family,   we  prove that the cardinal of the set of holomorphic gauge fields on any $B$-brane over the projective space is less than $2.$

\end{abstract}


%%%%%%%%%%%%%%%%%%%%%%%%%%%%%%%%%%%%%%%%%%%%%%%%%%%%%%%%%

\section {Introduction} \label{S:intro}
Given a holomorphic vector bundle $W$ over a complex manifold $Y,$  a
 connection on $W$ is holomorphic if the covariant derivative of any holomorphic section of $W$ is  also holomorphic.
Thus, the holomorphic connections are compatible with the holomorphic structures.

 Sixty-seven years ago, Atiyah  initiated the study of these connections in this 
context; in the category of holomorphic vector bundles   \cite{Atiyah}. Our purpose is to extend this concept to objects  of more general categories.

 But to which categories?
 The framework of vector bundles has some homological shortcomings. The category $\mathbf{Vec}(Y)$ of holomorphic vector bundles over the complex manifold $Y$  
 %of vector bundles on  the manifold  $Y$
  is not abelian: not every morphism has  cokernel. In fact, the cokernel of a morphism of vector bundles is a sheaf.
 
 A natural generalization would be to move to the category of sheaves.
% ; i.e the category of ${\mathscr O}_Y$-modules $\mathbf{Mod}_{{\mathscr O}_Y}$. 
However, there are sheaves so ``bad'' that they are even be supported on Cantor sets. 
%But a general sheaf can possess ``wild singularities''.
It is therefore advisable to restrict oneself to sheaves with``non-wild singularities''. The coherent sheaves are closely related with the geometry of the underlying space; furthermore,
 %It is therefore advisable to restrict oneself
% to coherent  sheaves.  
 the singularity locus of such a sheaf is a  subvariety with  codimension $\geq 1.$ The category of coherent sheaves over $Y$ will be denoted $\mathbf{Coh}(Y).$

Often when thinking of a single  sheaf, we are probably actually looking at a complex of sheaves. These complexes arise from injective resolutions, from
the differential forms, from the complexes of chains or cochains, etc. A complex 
of sheaves over $Y$ is a sequence of morphisms of sheaves
\begin{equation}\label{Sdot}
S^{\bf\centerdot}:\;\;\dots \to S^i\overset{d^i}{\to}S^{i+1}\to\cdots 
\end{equation}
satisfying $d^i\circ d^{i-1}=0.$
 The  bounded complexes of coherent sheaves over $Y$ are the objects of bounded derived category $D^b(Y).$

 On the other hand, from a mathematical point of view, a $B$-brane on   the    complex manifold $Y$ is an object of 
% the bounded derived category 
 $D^b(Y)$ \cite{Aspin, Aspin-et};  that is, a bounded complex of coherent sheaves. Thus, the   simplest $B$-branes are
 the objects of $\mathbf{Vec}(Y).$ We will extend the  concept of  holomorphic connection in the category $\mathbf{Vec}(Y)$ to the objects of $D^b(Y).$ In mathematical physics terms, we will define holomorphic gauge fields on $B$-branes.
 
  We will first provide a version of the definition of a holomorphic connection on a vector bundle that is suitable for extension to $B$-branes (Definition \ref{Def2}).
   We will then define the concept of a gauge field on $B$-branes (Definition \ref{Def:GaugeBrane}) in such a way that, when particularized to vector bundles, it coincides with the notion of a holomorphic connection.

   We will prove that the cardinal of the set of holomorphic gauge fields on any $B$-brane over ${\mathbb P}^n$ is $<2$ 
 (Theorem \ref{Thm1}). The proof of this theorem is based in two facts: 
 \begin{enumerate}
 \item The existence of a family $\{S_1,\dots, S_{n+1}\}$ of sheaves on the projective space ${\mathbb P}^n,$ which generates the category $D^b({\mathbb P}^n)$ and such that ${\rm Supp}\,S_{i}\cap{\rm Supp}\,S_j=\emptyset,$ for $i\ne j.$
 \item The vanishing of the Hodge cohomology groups $H^{1,0}({\mathbb P}^k).$
\end{enumerate}

A consequence of a celebrated theorem by Beilinson is the well-known fact that the   set  $\{{\mathscr O}_{{\mathbb P}^n}(-a);\; a=0,\dots, n\}$ generates $D^b({\mathbb P}^n)$ \cite[Cor. 8.29]{Huybrechts}. As the supports of these generators  are not disjoint, this set of generators is not suitable to prove our Theorem \ref{Thm1}.

This article is organized as follows.
%In Section \ref{S:Derived}
 The family of generators $\{S_i\}$ is constructed in Section \ref{S:Derived}. In this section, we briefly review some concepts from derived category theory that are necessary for the construction of 
% the aforementioned 
 that set of generators. The omitted details can be found in classical references such as \cite{Ge-Ma, Kas-Sch, Weibel}.

 In Section \ref{S:Holomorphic}, we define the holomorphic gauge fields on $B$-branes and prove Theorem \ref{Thm1}.

%%%%%%%%%%%%%%%%%%%%%%%%%%%%%%%%%%%%%%%%%%%%%%%%%%%%%%%%%%%%%%%%%%%%

%%%%%%%%%%%%%%%%%%%%%%%%%%%%%%%%%%%%%%%%%%%%%%%%%%%%%%%%%%%%%%%%%%%%%%%%%%%%%%%%%%%%%%%%%%%%%%%%%%%%%%%%%%%%%%%%%%%%%%%%%%%%%%%%%%%%%%%%
%%%%%%%%%%%%%%%%%%%%%%%%%%%%%%%%%%%%%%%%%%%%%%%%%%%%%%%%%%%%%%%%%%%%%%%%%%%%%%%%%%%%%%%%%%%%%%%%%%%%%%%%%%%%%%%%%%%%%%%%%%%%%%%%%%%

\section{Family of generators of $D^b({\mathbb P}^n)$}\label{S:Derived}
A morphism of complexes 
$(B^{\centerdot},\,d_B)\to (C^{\centerdot},\,d_C)$ in an abelian category ${\mathbf A}$ is a quasi-isomorphism, if it induces isomorphisms between the cohomologies $H^i(B^{\centerdot})\to H^i(C^{\centerdot}).$

The derived category $D^b({\mathbf A})$ has as objects the bounded complexes of ${\mathbf A}.$
The morphisms in $D^b({\mathbf A})$ are the morphisms of complexes in ${\mathbf A}$ together with the inverses of the quasi-isomorphisms \cite{Gefand-Manin2, Ge-Ma}. Thus,  
  every quasi-isomorphism     becomes an isomorphism  in the derived category.

\smallskip

%\noindent
\paragraph{\it The functor ${\mathbf A}\to D^b({\mathbf A})$}\;
Given $E\in{\rm Ob}({\mathbf A}),$ one defines the complex $Q(E)\in{\rm Ob}(D^b({\mathbf A}))$  by
$Q(E)^0=E$ and $Q(E)^p=0,$ for all $p\ne 0.$ The functor $Q:{\mathbf A}\to
D^b({\mathbf A})$ is fully faithfull \cite[p. 164]{Ge-Ma}:
\begin{equation}\label{fully-faithfull}
{\rm Hom}_{\mathbf A}(E_1,\,E_2)\simeq{\rm Hom}_{D^b({\mathbf A})}(Q(E_1),\,Q(E_2)).
\end{equation}

\smallskip

%%%%%%%%%%%%%%%%%%%%%%%%%%%%%%%%%%%%%%%%%%%%%%%%%%%%%%%%%%%%%%%%%%%%%%%%%
%%%%%%%%%%%%%%%%%%%%%%%%%%%%%%%%%%%%%%%%%%%%%%%%%%%%%%%%%%%%%%%%%%%%%%

%\noindent
\paragraph{\it The complex ${\rm Hom}$} \;
Given the complexes $B^{\centerdot}$ and $C^{\centerdot}$ in the category ${\mathbf A},$ one defines the  Hom complex ${\rm Hom}^{\centerdot}(B^{\centerdot},\,C^{\centerdot})$ by   (see \cite[page 17]{Iversen}) 
\begin{equation}\label{HomComplex}
{\rm Hom}^{m}(B^{\centerdot},\,C^{\centerdot})=\prod_{i\in{\mathbb Z}}{\rm Hom}_{\mathbf A}\big(B^i,\,C^{i+m}\big)
\end{equation}
with the differential  $d_H$
\begin{equation}\label{deltamg1}
(d_H^mg)^p=d_C^{m+p}g^p+(-1)^{m+1}g^{p+1}d_B^p.
\end{equation}
Furthermore, if   ${\rm Hom}^{\centerdot}(B^{\centerdot},\,C^{\centerdot})$ is the complex $0,$ then 
\begin{equation}\label{fully}
{\rm Hom}_{D^b({\mathbf A})}(B^{\centerdot},\,C^{\centerdot})= 0.
\end{equation}
%{\rm Hom}^{\centerdot}(B^{\centerdot},\,C^{\centerdot}).$

\smallskip

There are other two important constructions that can be carried out with the complexes in ${\mathbf A}$: The shifting and the mapping cone.

\smallskip

 \paragraph{\it The shifting}\; Given the complex $(A^{\centerdot},\,d_A)$ in the category ${\mathbf A}$   and $k\in{\mathbb Z},$ we denote  by $A[k]^{\centerdot}$ the complex $A^{\centerdot}$ shifted $k$ on the left; i. e.
$$A[k]^i=A^{k+i},\;\; d^i_{A[k]}=(-1)^kd^{i+k}_A.$$

\smallskip

 \paragraph{\it The mapping cone}\; 
  With the morphism  
  $h : A^{\centerdot} \to B^{\centerdot}$ of complexes
 one can define a new complex ${\rm Con}(h)^{\centerdot}$, called the mapping cone of $h,$ as follows:
$${\rm Con}(h)^{i}=A^{i+1}\oplus B^i,\;\;\;
d^i_{{\rm Con}} =\begin{pmatrix} -d^{i+1}_A & 0\\
h^{i+1} & d^i_B
\end{pmatrix}
$$

%Thus, %with the natural morphisms of complexes, 
 With the inclusion and the projection, we can form the following sequence of complexes.
\begin{equation}\label{triangle}
A^{\centerdot}\overset{h}{\longrightarrow}B^{\centerdot}\overset{i(h)}{\longrightarrow} {\rm Con}(h)^{\centerdot}\to A[1]^{\centerdot}\overset{-h}{\longrightarrow} B[1]^{\centerdot},
\end{equation}
where in the third position is the complex ${\rm Con}$ of the first morphism.  The sequence (\ref{triangle}) is called a {\it distinguished triangle} and sometimes is written 
 $$A^{\centerdot}\overset{h}{\longrightarrow} B^{\centerdot}\to {\rm Con}(h)^{\centerdot}\overset{+1}{\to}.$$
 
Repeating the construction with the morphism  $i(h),$  we obtain the corresponding distinguished triangle
\begin{equation}\label{triangle1}
B^{\centerdot}\overset{i(h)}\longrightarrow{}{\rm Con}(h)^{\centerdot}\to {\rm Con}(i(h))^{\centerdot}\to B[1]^{\centerdot}\to .
\end{equation}
It can be proved that the complex $A[1]^{\centerdot}$ is quasi-isomorphic to ${\rm Con}(i(h))^{\centerdot}$ \cite[Lem. 1.4.2]{Kas-Sch}.
 As quasi-isomorphisms in ${\mathbf A}$ become isomorphisms in $D^b({\mathbf A}),$  one can identify $A[1]^{\centerdot}$ and ${\rm Con}(i(h))^{\centerdot}$ in the derived category. Hence, in addition to the distinguished triangle (\ref{triangle}), we have the following  distinguished triangle in $D^b({\mathbf A}).$
 \begin{equation}\label{triangle2}
B^{\centerdot}\overset{i(h)}\longrightarrow{}{\rm Con}(h)^{\centerdot}\to A[1]^{\centerdot}\to B[1]^{\centerdot}\to .
\end{equation}
 That is, one has the folowing proposition:
\begin{Prop}\label{translation}
If the  sequence  $A^{\centerdot}\to  B^{\centerdot}\to C^{\centerdot}\to A[1]^{\centerdot}$ is  distinguished triangle in $D^b({\mathbf A}),$ then   so are the following
 $$ B^{\centerdot}\to  C^{\centerdot}\to A[1]^{\centerdot} \to B[1]^{\centerdot}\;\;  \text{and}\;\;
  C^{\centerdot}\to A[1]^{\centerdot}\to B[1]^{\centerdot}\to C[1]^{\centerdot}.$$   
\end{Prop}

 \smallskip
 
  Given an exact sequence of complexes in  ${\mathbf A}$ %of  coherent sheaves
  \begin{equation}\label{WCone}
  0\to A^{\centerdot}\overset{h}\longrightarrow B^{\centerdot} \to C^{\centerdot}\to 0,
  \end{equation}
  one can prove that the complexes ${\rm Con}(h)^{\centerdot}$ and $C^{\centerdot}$ are quasi-isomorphic \cite[Prop. 1.7.5]{Kas-Sch}.  Thus, they are isomorphic in the derived category $D^b({\mathbf A}),$  and  we have the following  proposition:
  \begin{Prop}\label{Exact-Cone}
  The short exact sequence 
   $0\to A^{\centerdot} \longrightarrow B^{\centerdot} \to C^{\centerdot}\to 0$ of complexes in ${\mathbf A}$ defines the following distinguished triangle
   % the following   distinguished triangle 
 $A^{\centerdot}\to B^{\centerdot}\to C^{\centerdot} \overset{+1}{\to}$
 in   $D^b({\mathbf A}).$
 \end{Prop}  
 
 %is a distinguished triangle as well. 
 %As the third term of a distinguished triangle is a cone map,
  %$B^{\centerdot}[1]$ is the cone map of a morphism
%$C^{\centerdot}\to A^{\centerdot}[1].$ That is, in the category %$D^b({\mathbf A})$, 
%the objects $B^{\centerdot}[1]$ and $C^{\centerdot}[1]\oplus A^{\centerdot}[1]$ are isomorphic. That is,
%$ \begin{equation}\label{T=S+W}
% B^{\centerdot} \simeq C^{\centerdot} \oplus A^{\centerdot}.\end{equation}

\smallskip

In the following, the category ${\mathbf A}$ will   often  be
 that of coherent sheaves on the manifold $Y,$  in which case the corresponding derived category  will be denoted $D^b(Y).$ 
 %or that of finitely generated $R$-modules ($R$ being a commutative ring).

\subsection{Generators of a derived category}
 In  the derived category there are two fundamental operations built in: \begin{itemize}
 \item The shifting of the complexes; that is the construction of the complexes ${A}[k]^{\centerdot}$ from ${A}^{\centerdot}.$
 \item The mapping cone operation; that is, the construction from a morphism $h:A^{\centerdot}\to B^{\centerdot}$ of the 
 ${\rm Con}(h)^{\centerdot}=A[1]^{\centerdot}\oplus B^{\centerdot}.$
 \end{itemize}

 Given a set ${\sf G}$ of objects of $D^b(\mathbf{A})$, the category $\mathbf{Cat}({\sf G})$ generated by ${\sf G}$ is the smallest full subcategory of $D^b(\mathbf{A}),$ such that
 \begin{enumerate}
 \item  It contains  ${\sf G}.$ %\subset{\mathbf Cat}({\sf G}).$
 \item  It is closed under shifting.
 \item It is closed under the mapping cone construction. That is,
 if $ A^{\centerdot} \to B^{\centerdot} \to  C^{\centerdot}  \overset {+1}{\to}$   is  dist.  triangle and  $A^{\centerdot}, 
 B^{\centerdot}\in \mathbf {Cat}({\sf G}),$   then  $C^{\centerdot}\in \mathbf{Cat}({\sf G}).$  
 \end{enumerate}
 
 From Proposition \ref{translation}, it follows the proposition:
 \begin{Prop}\label{Three}
 Each term of a dist. triangle $R^{\centerdot}\to S^{\centerdot}\to T^{\centerdot}\overset{+1}{\to}$ belongs to the category generated by the other two terms. 
 \end{Prop}

 As a consequence of Proposition \ref{Exact-Cone}, one has:
 % According to (\ref{exact-dist}), we can state the following proposition
 \begin{Prop}\label{Prop:Exact-Dist} 
   Given   the exact sequence    (\ref{WCone}),   then the complex $B^{\centerdot}$ belongs to the subcategory of $D^b({\mathbf A})$ generated by   $A^{\centerdot}$ and $C^{\centerdot}.$ 
 \end{Prop}
 
 \smallskip
 It is not difficult to prove that: If ${\sf G}$  is a set of objects of  $\mathbf{A}$ that generates $D^b({\mathbf A}),$  then any object of $D^b({\mathbf A})$    is isomorphic to a complex $F^{\centerdot},$ 
with
 \begin{equation}\label{formF}
 F^p=\bigoplus_{j}G_{pj}
 \end{equation}
 a finite direct sum of elements of ${\sf G}.$

 \smallskip
 
 \subsubsection{Affine varieties}

Given     a commutative ring $R$. We denote by $\mathbf{Mod}_R$ the category of finitely generated $R$-modules.
Let ${\sf G}=\{R\}$ be the singleton consisting of the $R$-module $R.$ The finite direct sum of copies of $R$ is an object of $\mathbf{Cat}({\sf G}),$ obviously.
   
 Given an object $M$ of $\mathbf{Mod}_R,$   according with the Hilbert's syzygy theorem, there exists a free resolution  
$$0\to F_k \overset{h_k}{\longrightarrow} F_{k-1} \overset{h_{k-1}}{\longrightarrow}\cdots \to F_1\overset{h_1}{\longrightarrow} F_0\to M\to 0,$$
where the $F_i$ are finite direct sums of $R$'s. 
 
From Propositions \ref{Exact-Cone} and \ref{Three}  applied to the short exact sequence
$$0\to F_k\to F_{k-1}\to {\rm Im}(h_{k-1})\to 0,$$ 
we deduce that, ${\rm Im}(h_{k-1})={\rm Ker}(h_{k-2})$ belongs to $\mathbf{Cat}({\sf G}).$ 

The induction  applied to the short exact sequences 
$$0\to {\rm Im} (h_{i+1})\to F_i\to{\rm Ker}(h_{i-1})\to 0$$
proves that $M\in \mathbf{Cat}({\sf G}).$   Thus,    $R$ is a generator of the derived category $D^b(\mathbf{Mod}_R).$ 
%of finite generated $R$-modules.

\smallskip

When $V$ is an affine variety, it is well-known that   the  category $\mathbf{Coh}(V)$, of coherent ${\mathscr O}_V$-modules, is equivalent to the category $\mathbf{Mod}_R$ of finitely generated $R$-modules,    $R$ being the ring of global sections of the structure sheaf ${\mathscr O}_V$ \cite[Ch. III, Sect 1]{Mumford}.
The equivalence is defined by the functor ``global section''
% $$\Gamma(U,\,-\): 
$${ H}\in \mathbf{Coh}(V)\mapsto \Gamma(V,\,H)\in\mathbf{Mod}_R.$$

 As a consequence of this equivalence, it follows the following proposition:
 \begin{Prop}\label{Prop:derivU}
If $V$ is an affine variety, then the derived category $D^b(V),$ of coherent sheaves on $V,$ is generated by the  sheaf ${\mathscr O}_V.$
\end{Prop}

%\smallskip

\subsubsection{Hypersufaces with affine complement}
Let $N$ be a hypersurface of the complex manifold $Y,$ $N\overset{i}{\hookrightarrow} Y.$ We set $V$ for the open 
$V:=Y\setminus N\overset{j}{\hookrightarrow} Y.$  We have the corresponding direct image functors and its adjoints
$$\mathbf{Mod}_{{\mathscr O}_N}\overset{i_*}{\longrightarrow} \mathbf{Mod}_{{\mathscr O}_Y}\overset{j_!}{\longleftarrow}
\mathbf{Mod}_{{\mathscr O}_V},\;\;\; \mathbf{Mod}_{{\mathscr O}_N}\overset{i^*}{\longleftarrow} \mathbf{Mod}_{{\mathscr O}_Y}\overset{j^!}{\longrightarrow}
\mathbf{Mod}_{{\mathscr O}_V}  $$
$i^*$ is the left adjoint of the direct image functor $i_*,$ and $j^!$ is the right adjoint of the proper direct image functor $j_!.$

One can consider the following three endofunctors of the category  $\mathbf{Mod}_{{\mathscr O}_Y},$
$$j_!j^!,\;i_*i^*,\;{\rm id}:\mathbf{Mod}_{{\mathscr O}_Y}  \longrightarrow
\mathbf{Mod}_{{\mathscr O}_Y}.$$
The adjuntion relations determine
 the following  natural transformations between these functors 
 $$j_!j^!\Longrightarrow {\rm id}\Longrightarrow i_{*}i^*,$$ 
 which in turn give rise,
  for each $S\in \mathbf{Coh}(Y),$ to the exact sequence  \cite[page 110]{Iversen}
\begin{equation}\label{exactjj}
0\to j_!j^!S\rightarrow {S}\rightarrow i_*i^*S\to 0.
\end{equation}
Furthermore, $j^!S\in\mathbf{Coh}(V)$ and $i^*S\in \mathbf{Coh}(N).$
From Proposition \ref{Prop:Exact-Dist}, it follows that
$S$ belongs to $\mathbf{Cat}\{j_!j^!S,\,i_*i^*S\},$ the subcategory of $D^b(Y)$ generated by these two sheaves. As $S$ is an arbitrary object
of $ \mathbf{Coh}(Y),$   we have the following proposition:
\begin{Prop}\label{Prp:CohZU}
Let ${\sf G}_V$ (${\sf G}_N$) be  a set of generators of $ D^b(V)$ (respec. $D^b(N)$), then the set
$(j_!{\sf G}_V)\cup(i_*{\sf G}_N)$ generates $D^b(Y).$
\end{Prop}

%\medskip

\subsubsection{Aplication to $D^b({\mathbb P}^n)$}
In ${\mathbb P}^n=:N_0$ we consider the hypersurface
$N_1=\{[x_0:\dots:x_n]\in {\mathbb P}^n; \,x_0=0\}={\mathbb P}^{n-1}.$ We set $V_1$  for the affine subvariety  
$V_1:=N_0\setminus N_1\simeq {\mathbb C}^n.$ One has the inclusions
$$N_1\overset{i_1}{\hookrightarrow} N_0 \overset{j_1}{\hookleftarrow} V_1.$$
Denoting by ${\sf G}_{N_1}$ a set of generators of  $D^b(N_1)$, from Proposition \ref{Prp:CohZU} together with Proposition \ref{Prop:derivU}, it follows that a set of generators of $D^b(N_0)$ is 
${\sf G}_{N_0}:=\{ j_{1!}{\mathscr O}_{V_1} \}\cup i_{1*}{\sf G}_{N_1}.$

 Next, we consider $N_2=\{[0:0:x_2:\dots :x_n]\in N_1\}.$ We set $V_2:=N_1\setminus N_2$ and 
  $N_2\overset{i_2}{\hookrightarrow} N_1 \overset{j_2}{\hookleftarrow} V_2.$
 Then  ${\sf G}_{N_1}=\{ j_{2!}{\mathscr O}_{V_2} \}\cup i_{2*}{\sf G}_{N_2}$ generates $D^b(N_1),$ assumed that ${\sf G}_{N_2}$ is a set of generators of $D^b(N_2).$ Thus,
 \begin{equation}\label{sfGn0}
 {\sf G}_{N_0}=\{ j_{1!}{\mathscr O}_{V_1} \}\cup
 \{ i_{1*}j_{2!}{\mathscr O}_{V_2} \}\cup i_{1*}i_{2*}{\sf G}_{N_2}.
 \end{equation}
 
 This process can be repeated $n$ times. Then $N_n$ is the singleton consisting of the point 
 %reduces to the point
  $z:=[0:\dots :0:1].$ We have the following diagram of inclusions
$$\xymatrix{N_n\ar[r]^{i_n} & N_{n-1}\ar[r]^{i_{n-1}} & \cdots\ar[r]^{i_3}  & N_2\ar[r]^{i_2} & N_1\ar[r]^{i_1}& N_0
 \\
 {}& V_{n}\ar[u]^{j_n} & {}  & V_3\ar[u]^{j_3} & V_2\ar[u]^{j_2}  & V_1\ar[u]^{j_1} 
  }$$
 For $k=1,\dots, n$ we let $\iota_k=i_{1}\circ\dots\circ i_{k}$ and
 $\iota_0={\rm id}. $  In this way,  one has the following ${\mathscr O}_{{\mathbb P}^n}$-module 
 $S_k:=\iota_{k-1 *}(j_{k!}{\mathscr O}_{V_k}),$ $k=1,\dots,n.$ With this   notation, (\ref{sfGn0}) reads
 $${\sf G}_{{\mathbb P}^n}=\{S_1,\,S_2\}\cup \iota_{2*}{\sf G}_{N_2}.$$ 

 The trivial  derived  category $D^b(N_n)$   is generated by the stalk  ${\mathscr O}_{{\mathbb P}^n,z}$ of ${\mathscr O}_{\mathbb{P}^n}$ at $z$.       Let $S_{n+1}$    be the skyscraper sheaf on ${\mathbb P}^n$ at  the point  $z.$ That is, $S_{n+1}=\iota_{n_*}({\mathscr O}_{{\mathbb P}^n,z})$. Then we can state: 

 \begin{Prop}\label{genrators Pn}
The family $\{S_1,\dots,S_{n+1}\}$ is a set of generators of the derived category $D^b({\mathbb P}^n).$
 \end{Prop}
 
 \begin{Rems}\label{Rem1}
For $k=1,\dots,n$ the sheaf $S_k$ is supported in $V_k$ and $S_{n+1}$ is supported in the point $N_{n}.$ On the other hand, 
 $V_k\subset N_{k-1}$ and $V_i\cap N_{k-1}=\emptyset$ for $i\leq k-1.$ In particular, if $r< k$, then   $V_r\cap V_k=\emptyset.$
%N_r\setminus \cup_{i=k}^rV_i.$ 
Hence, if $k\ne k'$ then 
\begin{equation}\label{disjoint}
{\rm Supp}\, S_k\cap {\rm Supp}\, S_{k'}=\emptyset.
\end{equation}
 \end{Rems}
 
 %We will prove the following Proposition:
By $\Omega^1$ we denote the sheaf of holomorphic $1$-forms on ${\mathbb P}^n.$ For the sake of simplicity, we set  ${\mathscr O}$ for ${\mathscr O}_{{\mathbb P}^n}$  and   $\Omega^1(S_i):=\Omega^1\otimes_{{\mathscr O}}S_i.$
\begin{Prop}\label{Lem1}
 For all $i,j\in\{ 1,\dots,n+1\}.$
$${\rm Hom}_{{\mathscr O}}(S_i,\,\Omega^1(S_j))=0,$$   
\end{Prop}
{\it Proof.} 
For $i\ne j$, by (\ref{disjoint}) 
%the above Remark
$${\rm Hom}_{\mathscr O}(S_i,\, \Omega^1 (S_j))=0.$$
On the other hand, $S_1=j_{1!}{\mathscr O}_{V_1}$; that is, it is isomorphic to the invertible 
%locally free 
${\mathscr O}$-module ${\mathscr O}(-N_1)$. Thus, tensoring by the dual $S_1^{\vee}$ of $S_1,$ one obtains
 $${\rm Hom}_{\mathscr O}(S_1,\, \Omega^1\otimes_{\mathscr O} S_1)\simeq
 {\rm Hom}_{\mathscr O}({\mathscr O},\, \Omega^1)\simeq H^0({\mathbb P}^n,\,\Omega^1).$$
This cohomology group vanishes \cite[page 4]{O-S-S}. Similarly, one can prove that
${\rm Hom}_{\mathscr O}(S_k,\, \Omega^1(S_k))=0,$ for all $k$.\qed

%%%%%%%%%%%%%%%%%%%%%%%%%%%%%%%%%%%%%%%%%%%%%%%%%%%%%%%%%%%%%%%%%%%%%%%%%%%%%%%%%%%%%%%%%%%%%%%%%%%%%%%%%%%%%%%%%%%%%%%%%%%%%%%%%%%%%%%%%%%%%%%%%%%%%%%%

\section{Holomorphic gauge fields}\label{S:Holomorphic}

\subsection{Holomorphic connections on a vector bundle}
As we mentioned in the Introduction, our purpose is to define holomorphic gauge fields on $B$-branes, extending the concept of   holomorphic connection on vector bundles.

First of all, let us recall the definition of a holomorphic connection on a holomorphic vector bundle $W\to Y.$ We denote also 
 by $W$ the locally free sheaf consisting of the sections of the 
 %holomorphic 
 vector bundle $W$. By $\Omega^1_Y$ is denoted the sheaf of holomorphic $1$-forms on $Y$.
  A holomorphic connection on $W$   is a morphism of abelian sheaves 
 $$\nabla:W\to \Omega^1(W)=\Omega^1_Y\otimes_{{\mathscr O}_Y}W,$$ 
 satisfying
 $$\nabla(f\sigma)=\partial f\otimes\sigma+f\sigma,$$
 where $f\in {\mathscr O}_Y$ is a function of the structure sheaf of $Y$ and $\sigma $ a holomorphic section of $W.$

 This definition admits another equivalent formulation, more appropriate for extension to $B$-branes, by means the $1$-jet bundle $J^1(W)$ of $W.$ The $1$-jet bundle is the abelian sheaf
 $$J^1(W)=W\oplus\Omega^1(W)$$
 endowed with the following ${\mathscr O}_Y$-module structure
 $$f\cdot(\sigma\oplus\beta)=f\sigma\oplus(\partial f\otimes\sigma+f\beta).$$
Denoting by $p:J^1(W)\to W$ the projection, one has the following exact sequence of ${\mathscr O}_Y$-modules
\begin{equation}\label{exactsequenceJW}
0\to\Omega^1(W)\to J^1(W)\overset{p}{\to} W\to 0.
\end{equation}

On the other hand, the inclusion 
$$t:\sigma\in W\mapsto \sigma \oplus 0\in J^1(W)$$
is a morphism of {\it abelian} sheaves such that $p\circ t={\rm id}.$

Given $\varphi\in {\rm Hom}_{{\mathscr O}_Y}(V,\,J^1(V))$ a right inverse of $p;$  that is, such that $p\circ \varphi={\rm id}.$  Then $p\circ (\varphi-t)=0$ and thus $\varphi- t$ factors uniquely through  ${\rm Ker}(p) =\Omega^1(W),$ defining the morphism $\nabla$ that is a connection on $W$.

$$ \xymatrix{0 \ar[r] & \Omega^1(W)\ar[r] & J^1(W)\ar[r]^p& W \ar[r] & 0 \\ 
   & & & W \ar[ul]_{\varphi-t} \ar@{-->}[ull]^{\nabla} &  } $$
   %A holomorphic connection is a right inverse of $p.$ In other words,
Hence, one can give a new definition of holomorphic connection equivalent to the preceding one.
\begin{Def}\label{Def2}  The holomorphic  connections on $W$ are the elements of the following set
\begin{equation}\label{Def:connection}
\{\varphi\in {\rm Hom}_{{\mathscr O}_Y}(W,\,J^1(W)); \; p\circ \varphi={\rm id}\}.
\end{equation}
\end{Def}
  That is, a holomorphic connection is a splitting of the  exact sequence (\ref{exactsequenceJW}).
  
The induced Ext exact sequence reads
\begin{align} \notag 
 0&\to {\rm Hom}_{{\mathscr O}_Y}(W,\,\Omega^1(W))  \to {\rm Hom}_{{\mathscr O}_Y}(W,\,J^1(W))\to \\
&\to{\rm Hom}_{{\mathscr O}_Y}(W,\,W)\to {\rm Ext}^1_{{\mathscr O}_Y}(W,\,\Omega^1(W))\to \cdots \notag  
\end{align}

 The image of ${\rm id}\in{\rm Hom}_{{\mathscr O}_Y}(W,\,W)$ in ${\rm Ext}^1_{{\mathscr O}_Y}(W,\,\Omega^1(W))$ is the Atiyah class of $W.$ The second non trivial morphism in the Ext sequence is the map $\varphi\mapsto   p\circ\varphi.$ Thus, $p$ admits a right inverse iff the Atiyah class of $W$ vanishes. Furthermore, if there exist  two right inverses $\varphi$ and $\varphi'$ of $p$, then $p(\varphi-\varphi')=0.$ That is, 
 $$\varphi-\varphi'\in {\rm Hom}_{{\mathscr O}_Y}(W,\,\Omega^1(W)).$$
 One has the following proposition:
 \begin{Prop}\label{P:Atiyah}
 The holomorphic vector bundle $W$ admits a holomorphic connection iff its Atiyah class vanishes. When the set of holomorphic connections is non empty, it is an affine space associated to the finite dimensional vector space
 ${\rm Hom}_{{\mathscr O}_Y}(W,\,\Omega^1(W)).$
 \end{Prop}

%\medskip

%%%%%%%%%%%%%%%%%%%%%%%%%%%%%%%%%%%%%%%%%%%%%%%%%%%%%%%%%%%%%%%%%%%%%%%%%%%%

\subsection{Holomorphic gauge fields on a $B$-brane.}

A bounded   complex $S^{\centerdot}$ 
%(\ref{Sdot}) 
of coherent sheaves over $Y$ is an object in the derived category $D^b(Y)=D^b(\mathbf{Coh}(Y)).$ One also has the corresponding $1$-jet complex 
$$J^1(S^{\centerdot})=S^{\centerdot}\oplus \Omega^1(S^{\centerdot})\overset{p^{\centerdot}}{\longrightarrow} S^{\centerdot}.$$
 % morphisms of complexes
%$ J^1(S^{\centerdot})\overset{\pi}{\to} S^{\centerdot}.$
According to Definition \ref{Def2},  
%\ref{Def:connection}, 
it is it is reasonable to define the gauge fields on the $B$-brane $S^{\centerdot}$  as the right inverses of the morphism $p^{\centerdot}$. That is,
\begin{Def}\label{Def:GaugeBrane} 
The holomorphic gauge fields on the $B$-brane $S^{\centerdot}$ are the elements of 
\begin{equation}\label{DefGauge}
\big\{\phi\in{\rm Hom}_{D^b(Y)}(S^{\centerdot},\,J^1(S^{\centerdot})); \; p^{\centerdot}\circ\phi={\rm id}   \big\}.
\end{equation}
\end{Def}
By (\ref{fully-faithfull}), this definition, when applied to complexes consisting of only one nontrivial term which is a locally free sheaf, coincides with Definition \ref{Def2}.
% for holomorphic vector bundles.
%   given in the previous section.

Although $D^b(Y)$ is not an abelian category the exact sequence of complexes of ${\mathscr O}_Y$-modules 
$$0\to\Omega^1(S^{\centerdot})\to J^1(S^{\centerdot}) \to  S^{\centerdot}\to 0,$$
gives rise, according to Proposition \ref{Exact-Cone}, to the following distinguished triangle in $D^b(Y)$
$$\Omega^1(S^{\centerdot})\to J^1(S^{\centerdot})  \to  S^{\centerdot}\overset{+1} {\longrightarrow}.
$$

Since ${\rm Hom}_{D^b(Y)}(S^{\centerdot},\,-\,)$ is a cohomological functor  \cite[Prop. 1.5.3]{ Kas-Sch}, from the above distinguished triangle we deduce the following exact sequence
 \begin{align}
 \notag 0 & \to  {\rm Hom}_{D^b(Y)}(S^{\centerdot},\, \Omega^1(S^{\centerdot}))\to {\rm Hom}_{D^b(Y)}(S^{\centerdot},\, J^1(S^{\centerdot}))\to \\ \notag
&\to {\rm Hom}_{D^b(Y)}(S^{\centerdot},\,  S^{\centerdot})\to
 {\rm Ext}^1(S^{\centerdot},\, \Omega^1(S^{\centerdot}))\to
 \end{align}
 The image of ${\rm id}\in {\rm Hom}_{D^b(Y)}(S^{\centerdot},\,  S^{\centerdot})$ in the space ${\rm Ext}^1(S^{\centerdot},\, \Omega^1(S^{\centerdot}))$ is ${\rm At}(S^{\centerdot})$ is the Atiyah class of the 
 $B$-brane $S^{\centerdot}$. The vanishing of ${\rm At}(S^{\centerdot})$ is a necessary and sufficient condition for the existence of holomorphic gauge fields on $S^{\centerdot}.$
From the exactness of the  Ext sequence, it also follows  the following proposition:
\begin{Prop}\label{Propafine}
If the set of holomorphic gauge fields on the $B$-brane ${S^{\centerdot}}$ is non empty, then it is an affine space associated to the finite dimensional vector space ${\rm Hom}_{D^b(Y)}(S^{\centerdot},\, \Omega^1(S^{\centerdot})).$
\end{Prop}

\subsection{Holomorphic gauge fields  on $B$-branes over ${\mathbb P}^n$}
By Proposition \ref{genrators Pn} and according to (\ref{formF}), any  object of $D^b({\mathbb P}^n)$ is isomorphic to one complex
 $F^{\centerdot}$ of the form $F^{p}=\bigoplus_jS_{pj},$ with $S_{pj}\in\{S_1,\dots, S_{n+1}  \}.$ Thus,
%On the other hand,
$\Omega^1(F)^{q}=\bigoplus_k\Omega^1(S_{qk}).$ 

From Proposition \ref{Lem1}, it follows that 
${\rm Hom}_{{\mathscr O}_{{\mathbb P}^n}}(F^p,\,\Omega^1(F^q))=0.$
 Hence, the complex 
${\rm Hom}^{\centerdot}(F^{\centerdot},\,\Omega^1(F)^{\centerdot})=0.$ By (\ref{fully}), we deduce 
%From (\ref{HomComplex}), one deduces
$${\rm Hom}_{D^b({\mathbb P}^n)}\big(F^{\centerdot},\, \Omega^1(F)^{\centerdot}\big)=0.$$
  The following theorem is a consequence of Proposition \ref{Propafine}.
\begin{Thm}\label{Thm1}
The number of holomorphic  gauge fields on an arbitrary $B$-brane over ${\mathbb P}^n$ is either  
$0$ or $1.$   
\end{Thm}

  Therefore, for a $B$-brane  over ${\mathbb P}^n$, either it is not possible to define a holomorphic covariant derivative of its sections, or only a single one can be defined. In particular,  on the sheaf ${\mathscr O}_{{\mathbb P}^n}$ the operator $\partial:{\mathscr O}_{{\mathbb P}^n}\to \Omega^1_{{\mathbb P}^n}$ is the only holomorphic connection on this $B$-brane.
  
  \medskip
  \begin{Rems}
  Let $Y$ be a complex manifold such that, there exists a ``tower'' 
  $Y=N_0\supset N_1\supset N_2\supset\dots \supset\{ {\rm point}\}$ of submanifolds of $Y$ satisfying:
  \begin{enumerate}
  \item $N_i$ is a divisor of $N_{i-1}$
  \item $N_{i-1}\setminus N_i$ is an affine variety
  \item The Hodge cohomology groups $H^{1,0}(N_i)=0$
  \end{enumerate}
  Then, mimicking the development made for the case of ${\mathbb P}^n,$ one can prove that the set of holomorphic gauge fields on any $B$-brane over $Y$
  has cardinal $<2.$ 
  \end{Rems}

\smallskip
%%%%%%%%%%%%%%%%%%%%%%%%%%%%%%%%%%%%%%%%%%%%%%%%%%%%%%%%%%%%%%%%%%%%%%%%%%%%%%%%%%%%%%%%%%%%%%%%%%%%%%%%%%%%%%%%%%%%%%%%%%%%%%%%%
%%%%%%%%%%%%%%%%%%%%%%%%%%%%%%%%%%%%%%%%%%%%%%%%%%%%%%%%%%%%%%%%%%%%%%%%%%%%%%%%%%%%%%%%%%%%%%%%%%%%%%%%%%%%%%%%%%%%%%%%%%%%%%%%%%

%%%%%%%%%%%%%%%%%%%%%%%%%%%%%%%%%%%%%%%%%%%%%%%%%%%%%%%%%%%%%%%%%%%%%%%%%%%%%%%%%%%%%%%%%%%%%%%%%%%%%%%%%%%%%%%%%%%%%%%%%%%%%%%%%%%%%%%%%

%%%%%%%%%%%%%%%%%%%%%%%%%%%%%%%%%%%%%%%%%%%%%%%%%%%%%%%%%%%%%%%%%%%%%%%%%%%%%%%%%%%%%%%%%%%%%%%%%%%%%%%%%%%%%%%%%%%%%%%%%%%%%%%%%%%%%%
%%%%%%%%%%%%%%%%%%%%%%%%%%%%%%%%%%%%%%%%%%%%%%%%%%%%%%%%%%%%%%%%%%%%%%%%%%%%%%%%%%%%%%%%%%%%%%%%%%%%%%%%%%%%%%%%%%%%%%%%%%%%%%%%%%%%%%%%%

\end{document}